
\input amstex
\documentstyle{amsppt}
\refstyle{A}
\magnification=\magstephalf
\nologo
\NoBlackBoxes
\hsize=6in

\font\elevenbx=cmbx10 scaled \magstep1


\topmatter
\title
Batalin-Vilkovisky Lie algebra structure on the loop homology of complex Stiefel manifolds
\endtitle
\rightheadtext{BV Lie algebra structure on loop homology}
\author
Hirotaka Tamanoi
\endauthor
\affil
University of California, Santa Cruz
\endaffil
\address
Department of Mathematics, University of California Santa Cruz, 
Santa Cruz, CA 95064
\endaddress
\email
tamanoi\@math.ucsc.edu
\endemail
\abstract
We determine the Batalin-Vilkovisky Lie algebra structure for the
integral loop homology of special unitary groups and complex Stiefel
manifolds. It is shown to coincide with the Poisson algebra structure
associated to a certain odd symplectic form on a super vector space
for which loop homology is the super algebra of functions. Over
rationals, the loop homology of the above spaces splits into a tensor
product of simple BV algebras, and it is shown to contain the super
Lie algebra $\pi sp$.
\endabstract
\toc
\head
1.  Introduction \dotfill 1
\endhead
\head
2.  Intersection and Pontrjagin products in the
   homology of $\text{\rm SU}(n+1)$ \dotfill 4
\endhead
\head
3.  Homology $S^1$-action on the Hopf algebra $H_*(L\text{\rm
    SU}(n+1))$  
\dotfill 7
\endhead
\head
4.  BV Lie algebra structure on the loop homology
   $\Bbb H_*\left(L\text{\rm SU}(n+1)\right)$ \dotfill 10
\endhead
\head
5. BV-algebra structure on the loop homology of complex
   Stiefel manifolds \dotfill 14
\endhead
\head
{} References \dotfill 16
\endhead
\endtoc
\endtopmatter

\document

\head
{\elevenbx \S 1\quad Introduction}
\endhead

The homology of the free loop space $LM$ of an oriented closed
manifold $M$ was shown to admit a structure of a Batalin-Vilkovisky
(BV) algebra in which an associative algebra structure, a
Lie algebra structure and a circle action are amalgamated into a
single coherent algebraic structure \cite{1}. Namely, the Lie bracket
appears as the deviation of homology $S^1$-action from being a
derivation with respect to the loop product which defines the
associative algebra structure. The loop product structure in loop
homology has been calculated for some spaces including spheres and
complex projective spaces \cite{3}. However, the nature of the Lie
algebra structure on loop homology has not been well understood,
except that one of the defining identities of BV algebras says that
the Lie bracket acts as a derivation on the underlying
associative algebra.  The purpose of this paper is to describe the full
BV algebra structure on the loop homology of complex
Stiefel manifolds. For a brief description of a general BV algebra, see the
beginning of section 4. 

Since complex Stiefel manifolds are homogeneous spaces of special
unitary groups, the description of BV Lie algebra structure on the
loop homology of the Lie group $\text{\rm SU}(n+1)$ is fundamental for
our purpose, and the results for Stiefel manifolds follow from it. In
general, for a Lie group $G$, the associative algebra structure of its
loop homology is given by a tensor product $\Bbb H_*(LG)\cong \Bbb
H_*(G)\otimes H_*(\Omega G)$, where $\Bbb H_*(G)\cong H^{-*}(G)$ is
the intersection homology ring with degree shift so that the
fundamental class has degree $0$, and $\Omega G$ is the group of based
loops. Although the associative algebra structure is the obvious one,
the BV Lie algebra structure on it is highly nontrivial because the
circle action on the loop group $LG$ mixes the homology of $G$ and
the homology of $\Omega G$. For the case of $G=\text{\rm SU}(n+1)$, we
show that the BV Lie bracket coincides with the Poisson bracket
associated to an odd symplectic form on a super (vector) space of
dimension $n|n$ whose super algebra of functions is the loop homology
$\Bbb H_*\left(L\text{\rm SU}(n+1)\right)$. Similar description is
valid for the loop homology of all complex Stiefel manifolds.

To describe the result, the associative algebra structure of the
integral loop homology is given by
$$
\Bbb H_*\left(L\text{\rm SU}(n+1)\right)\cong
\Lambda_{\Bbb Z}(\alpha_3, \alpha_5, \cdots, \alpha_{2n+1})\otimes
\Bbb Z[e_2, e_4, \cdots, e_{2n}],
\tag 1-1 
$$
where $\Bbb H_*\left(\text{\rm SU}(n+1)\right)\cong 
\Lambda_{\Bbb Z}(\alpha_3, \alpha_5, \cdots, \alpha_{2n+1})$ and
$H_*\left(\Omega\text{\rm SU}(n+1)\right)\cong 
\Bbb Z[e_2, e_4, \cdots, e_{2n}]$ are also coalgebras whose generators
have degree $|\alpha_{2k+1}|=-2k-1$, $|e_{2k}|=2k$ for $1\le k\le
n$. For sequences $I=(2i_1+1, 2i_2+1, \cdots, 2i_r+1)$ with $1\le
i_1<i_2<\cdots <i_r\le n$, and $J=(j_1,j_2,\cdots, j_n)$ with
nonnegative integer entries, let
$\alpha_I=\alpha_{2i_1+1}\alpha_{2i_2+1}\cdots \alpha_{2i_r+1}$ and
$e^J=e_{2}^{j_1}e_4^{j_2}\cdots e_{2n}^{j_n}$. We set $|I|=r$, the
length of $I$. 

We define a differential operator $D_{2\ell}$
acting on the polynomial ring $\Bbb Z[e_2, e_4, \dots, e_{2n}]$ by 
$$
D_{2\ell}=\sum_{k=\ell}^ne_{2k-2\ell}\frac{\partial}{\partial e_{2k}}
=\frac{\partial}{\partial e_{2\ell}}
+e_2\frac{\partial}{\partial e_{2\ell+2}}+\cdots+
e_{2n-2\ell}\frac{\partial}{\partial e_{2n}}, \quad 1\le\ell\le n, 
\tag1-2
$$
where we set $e_0=1$.  These $n$ differential operators mutually
commute. On the exterior algebra $\Lambda_{\Bbb Z}(\alpha_3, \alpha_5,
\cdots, \alpha_{2n+1})$, we have {\it odd} derivations
$\partial/\partial\alpha_{2\ell+1}$ for $1\le\ell\le n$ in the sense
that when we move $\partial/\partial\alpha_{2\ell+1}$ beyond
$\alpha_{2k+1}$, we get a minus sign.

The circle action $\Delta: S^1\times LM @>>> LM$ induces the BV
operator $\pmb{\Delta}: \Bbb H_*(LM) @>>> \Bbb H_{*+1}(LM)$ given by
${\pmb\Delta}(z)=\bold{\Delta}_*([S^1]\otimes z)$ for $z\in \Bbb
H_*(LM)$. Recall that the BV Lie bracket $\{\cdot,\cdot\}$ and the BV
operator $\pmb{\Delta}$ are related by a formula
$$
\pmb{\Delta}(a\cdot b)=\pmb{\Delta}(a)\cdot b + (-1)^{|a|}a\cdot
\pmb{\Delta}(b) +(-1)^{|a|}\left\{a,b\right\},
\tag1-3
$$
for arbitrary two elements $a,b$ in a BV algebra. 

\proclaim{Theorem A} 
\rm{(1)} The BV operator $\pmb{\Delta}$ on the loop homology 
$\Bbb H_*\left(L\text{\rm SU}(n+1)\right)$ is given by 
$$
\pmb{\Delta}(\alpha_Ie^J)=
\sum_{\ell=1}^n\frac{\partial \alpha_I}{\partial \alpha_{2\ell+1}}
\cdot D_{2\ell}(e^J).
$$

\rom{(2)} The BV Lie bracket in the loop homology 
$\Bbb H_*\left(L\text{\rm SU}(n+1)\right)$ is given by 
$$
(-1)^{|I|}\left\{\alpha_Ie^J, \alpha_Ke^L\right\}
=\pmb{\Delta}(\alpha_Ie^L)\cdot \alpha_Ke^J
+(-1)^{|I|}\alpha_Ie^L\cdot\pmb{\Delta}(\alpha_Ke^J).
$$
\endproclaim

It is interesting that the formula for the BV Lie bracket is given by
a very clean derivation formula of $\pmb{\Delta}$, but with $e$-terms
switched. To prove Theorem A, the homological $S^1$-action is first
calculated in the Hopf algebra $H_*\left(L\text{\rm SU}(n+1)\right)$
with Pontrjagin product with respect to which $\Delta$ is a
derivation. Its coalgebra structure is used to determine $\Delta$ on
Pontrjagin ring generators. Then the result is translated to the loop
homology. 
 
To have a better understanding of the BV Lie bracket, we take a
point of view that the algebra in (1-1) is a super algebra of
functions on a super vector space in which $e_2, e_4, \dots, e_{2n}$
are even coordinate functions and $\alpha_3, \alpha _5, \dots,
\alpha_{2n+1}$ are odd coordinate functions. See \cite{4} for super
algebras and super geometry. In the Hopf algebra
$H_*\left(\Omega\text{\rm SU}(n+1)\right)\otimes \Bbb Q$ over
rationals, there exist primitive elements $h_{2\ell}$ such that
$h_{2\ell}=e_{2\ell}+(\text{decomposables})$ for $1\le\ell\le n$ so
that $H_*\left(\Omega\text{\rm SU}(n+1)\right)\otimes \Bbb Q
\cong\Bbb Q[h_2,h_4,\cdots, h_{2n}]$. We regard $h_{2\ell}$'s as new
even coordinate functions, and as such a differentiation
$\partial/\partial h_{2\ell}$ makes sense for functions in
$e_{2k}$'s. On the super vector space $\Bbb Q^{n|n}$, we consider an
odd symplectic form 
$$
\omega=\sum_{\ell=1}^nd\alpha_{2\ell+1}\wedge dh_{2\ell}.    
$$
We consider the associated Poisson bracket on functions 
$[\cdot\ ,\cdot]_\omega: \Bbb H_*^{\Bbb Q}\otimes\Bbb H_*^{\Bbb Q} @>>> 
\Bbb H_*^{\Bbb Q}$, where $\Bbb H_*^{\Bbb Q}=\Bbb H_*\otimes\Bbb Q$, 
given by  $[F,G]_{\omega}=X_FG$, where $X_F$ is the Hamiltonian vector field
associated to $F$ by $\iota_X\omega=dF$. Concretely, the Poisson
bracket is given by 
$$
\left[\alpha_Ie^J, \alpha_Ke^L\right]_{\omega}
=\sum_{\ell=1}^n\left(
(-1)^{|I|}\frac{\partial (\alpha_Ie^J)}{\partial \alpha_{2\ell+1}}\cdot
\frac{\partial(\alpha_Ke^L)}{\partial h_{2\ell}} +
\frac{\partial(\alpha_Ie^J)}{\partial h_{2\ell}}
\cdot\frac{\partial(\alpha_Ke^L)}{\partial \alpha_{2\ell+1}}\right).
$$
 
The next theorem explains the nature of the BV Lie bracket.  

\proclaim{Theorem B} The Poisson bracket is well defined over $\Bbb
H_*$, and the BV Lie bracket in the loop homology $\Bbb
H_*\left(L\text{\rm SU}(n+1)\right)$ coincides with the Poisson
bracket associated to the symplectic form $\omega$. Namely, 
$\left\{\alpha_Ie^J, \alpha_Ke^L\right\}
=\left[\alpha_Ie^J, \alpha_Ke^L\right]_{\omega}$ for $\alpha_Ie^J,
\alpha_Ke^L\in\Bbb H_*$.
\endproclaim

It is satisfying that symplectic even coordinate functions
$h_{2\ell}$'s are given exactly by primitive elements in the function
algebra $\Bbb H_*^{\Bbb Q}=\Bbb H_*\otimes\Bbb Q$.

Let $\Bbb H(\ell)=\Lambda_{\Bbb Z}(\alpha)\otimes\Bbb Z[h]$
with $|\alpha|=-2\ell-1$ and $|h|=2\ell$ be a BV algebra 
associated to an odd symplectic form $d\alpha\wedge dh$ and with a BV
operator given by $\pmb\Delta(h^k)=0$ and $\pmb\Delta(\alpha
h^k)=kh^{k-1}$ for $k\ge0$. Since the symplectic vector space $(\Bbb
Q^{n|n},\omega)$ is a direct sum $\bigoplus_{\ell=1}^n(\Bbb
Q^{1|1},d\alpha_{2\ell+1}\wedge dh_{2\ell})$, over rationals, our BV
algebra $\Bbb H_*\left(L\text{SU}(n+1)\right)$ admits a corresponding
splitting. Namely, 

\proclaim{Corollary C} Over rationals, the loop homology of $\text{\rm
SU}(n+1)$ splits into a following tensor product as BV algebras\rom{:}
$$
\Bbb H_*\left(L\text{\rm SU}(n+1)\right)\otimes\Bbb Q\cong
\bigotimes_{\ell=1}^n\Bbb H(\ell)\otimes\Bbb Q.
$$
\endproclaim
The above corollary is related to a fact that over rationals, the Lie
group $\text{SU}(n+1)$ splits as a product of odd dimensional spheres,
and the BV algebra $\Bbb H(\ell)$ is precisely the loop homology of
$S^{2\ell+1}$ in view of Theorem 5-2.  

In general, for a BV algebra $\Bbb H_*$, the Lie bracket is a
derivation in each variable. For example, 
$$
\{a,b\cdot c\}=\{a,b\}\cdot c +(-1)^{|b|(|a|+1)}b\cdot\{a,c\},
$$
for any three elements $a,b,c\in\Bbb H_*$. This means that there is a
Lie algebra homomorphism from $\Bbb H_*$ to the Lie algebra
$\text{Der}(\Bbb H_*)$ of all derivations on the associative algebra
$\Bbb H_*$. In our case of $\Bbb H_*=\Bbb H_*\left(L\text{\rm
SU}(n+1)\right)$, the Lie algebra $\text{Der}(\Bbb H_*)$ is a free
$\Bbb H_*$-module of rank $n|n$.  Theorem B suggests that the
description of the Lie algebra structure on the loop homology
simplifies when tensored with $\Bbb Q$ so that primitive elements are
available. For a derivation $X\in\text{Der}(\Bbb H^{\Bbb Q}_*)$, let $\Cal L_X$
be its Lie derivative given by $\Cal L_X=d\circ\iota_X+\iota_X\circ d$.

\proclaim{Corollary D} 
\rom{(1)} The BV Lie algebra $\Bbb H_*^{\Bbb Q}/\Bbb Q\cdot 1$ is
isomorphic to the super Lie algebra of Hamiltonian derivations given
by $\{X\in\text{\rm Der}(\Bbb H_*^{\Bbb Q})\mid \Cal L_X\omega=0\}$. 

\rom{(2)} Let $Q\subset\Bbb H_*^{\Bbb Q}$ be the subspace spanned by
quadratic elements in $h_{2\ell}$'s and $\alpha_{2k+1}$'s. Then $Q$ is
a super Lie subalgebra isomorphic to the super Lie algebra
$\pi sp(n|n)$. 
\endproclaim
Here the super Lie algebra $\pi sp(n|n)$ consists of those linear maps
on a vector space with basis
$\{h_2,\cdots,h_{2n},\alpha_3,\cdots,\alpha_{2n+1}\}$ preserving the
odd symplectic form $\omega$ \cite{5}.

For the complex Stiefel manifold
$V_{n+1-k}(\Bbb C^{n+1})=\text{\rm SU}(n+1)/\text{\rm SU}(k)$ of
$n+1-k$ orthonormal frames in $\Bbb C^{n+1}$,
entirely analogous statements hold in which
$$
\Bbb H_*\left(V_{n+1-k}(\Bbb C^{n+1})\right)
\cong \Lambda_{\Bbb Z}(\alpha_{2k+1}, \cdots, \alpha_{2n+1})\otimes
\Bbb Z[e_{2k},\cdots, e_{2n}],
$$
with some obvious modifications on the operators $D_{2\ell}$'s.  

In this article, all homology and cohomology groups have integer
coefficients unless otherwise stated. 

The organization of this paper is as follows. After reviewing the
relationship of intersection product and Pontrjagin product in the
homology of $\text{SU}(n+1)$ in section 2, we compute the homological
$S^1$ action in the Pontrjagin ring in section 3. Here the Hopf
algebra is effectively used to determine this action. Then in section
4, we describe the full BV algebra structure of the loop homology of
$\text{SU}(n+1)$ using the result of section 3. After deriving the
formula for the BV Lie bracket, we describe its meaning in the context
of the odd symplectic structure.  In section 5, we describe BV algebra
structure on the loop homology of complex Stiefel manifolds.

\bigskip

\head
{\elevenbx \S2 \quad Pontrjagin and intersection products in 
$H_*\left(\text{\rm SU}(n+1)\right)$}
\endhead
 
\bigskip

In this section, we review the homology and cohomology groups of
$\text{\rm SU}(n+1)$ as Hopf algebras, with a discussion on a relation
between Pontrjagin product and intersection product in
$H_*\left(\text{\rm SU}(n+1)\right)$. We start with the Pontrjagin
structure. A basic reference for the topology of Lie groups is
\cite{6}. 

Let $\Sigma\Bbb CP^n$ be the reduced suspension of $\Bbb
CP^n=\text{\rm U}(n+1)/\left(\text{\rm U}(1)\times\text{\rm
U}(n)\right)$. Let a map $\beta_n:\Sigma\Bbb CP^n @>>> \text{\rm
SU}(n+1)$ be given by $\beta_n(t, [A])=A\cdot R_t(1,n)\cdot
A^{-1}\cdot R_t(1,n)^{-1}$, where $A\in\text{\rm U}(n+1)$ and
$R_t(1,n)=e^{\pi i t}\oplus e^{-\pi i t}I_n$. It is easy to check that
this definition is independent of the choice of the representative $A$
of the equivalence class $[A]$. Let $\bar{\beta}_n : \Bbb CP^n @>>>
\Omega \text{\rm SU}(n+1)$ be the adjoint of $\beta_n$. The the
Pontrjagin ring structures on the homology of $\text{\rm SU}(n+1)$ and
of its based loop space $\Omega\text{\rm SU}(n+1)$ are given by the
following well known theorem.

\proclaim{Theorem 2-1} \rom{(1)} Let $x_{2k+1}={\beta_n}_*([\Sigma\Bbb
CP^k])$ for $1\le k\le n$. Then the Hopf algebra structure of the 
homology of $\text{\rm SU}(n+1)$ is given by 
$$
H_*\left(\text{\rm SU}(n+1)\right)\cong \Lambda_{\Bbb Z}(x_3, x_5,
\cdots, x_{2n+1}),
$$
and all generators are primitive. 

\rom{(2)}  Let $e_{2k}={\bar{\beta}_{n*}}([\Bbb CP^k])$ for $1\le k\le
n$. Then the Pontrjagin ring structure and the coalgebra structure of
the homology of $\Omega \text{\rm SU}(n+1)$ are given by 
$$
H_*\left(\Omega \text{\rm SU}(n+1)\right)
\cong\Bbb Z[e_2,e_4, \cdots, e_{2n}], \qquad 
\phi_*(e_{2k})=\sum_{i+j=k}e_{2i}\otimes e_{2j}.
$$
\endproclaim

We also regard homology elements $x_{2k+1}$ and $e_{2k}$ as maps
$x_{2k+1}:\Sigma\Bbb CP^k @>>> \text{\rm SU}(n+1)$ and 
$e_{2k}: \Bbb CP^k @>>> \Omega\text{\rm SU}(n+1)$, which are adjoint
to each other, for $1\le k\le n$. 

Let $U$ be the sequence $U=(3,5,\cdots, 2n+1)$. For any sequence
$I=(2i_1+1, 2i_2+1,\cdots, 2i_r+1)$ with $1\le i_1<\cdots<i_r\le n$,
we let $x_I=x_{2i_1+1}\cdot x_{2i_2+1}\cdots x_{2i_r+1}$.  We orient
$\text{\rm SU}(n+1)$ by the fundamental class $x_U$. Let $y_{2k+1}$ be
the cohomology element dual to $x_{2k+1}$ with respect to the monomial
basis $\{x_I\}$. It is straightforward to dualize the Hopf algebra
structure of $H_*\left(\text{\rm SU}(n+1)\right)$ to get the Hopf
algebra structure of the cohomology. 

\proclaim{Theorem 2-2} As a Hopf algebra we have 
$$
H^*\left(\text{\rm SU}(n+1)\right)\cong
\Lambda_{\Bbb Z}(y_3,y_5,\cdots, y_{2n+1}),
$$ 
where all the exterior generators are primitive and the Kronecker
pairing is given by
$$
\langle y_I, x_J\rangle=(-1)^{\frac{|I|(|I|-1)}2}\delta_{I,J}.
\tag2-1
$$
\endproclaim

The sign above is the result of the usual sign convention in which
switching two odd degree objects gives a minus sign. 

Next we discuss the intersection product in homology. For an oriented
closed manifold $M^d$, the intersection product in homology can be
defined as the Poincar\'e dual of the cup product in
cohomology. Namely, let $D: H_*(M) @>{\cong}>> H^{d-*}(M)$ be the
Poincar\'e duality isomorphism given by $D(\alpha)\cap [M]=\alpha$ for
$\alpha\in H_*(M)$. Then for $\alpha,\beta\in H_*(M)$, their
intersection product $\alpha\circ\beta$ is defined by
$D(\alpha\circ\beta)=D(\alpha)\cup D(\beta)$. It is immediate that
$\beta\circ\alpha=(-1)^{(d-|\alpha|)(d-|\beta|)}\alpha\circ\beta$.

\definition{Remark} If $\alpha,\beta$ are homology classes represented by
submanifolds $K,N$ intersecting transversally, the orientation of
$\alpha\circ\beta=[K\cap N]$ is given as follows. At $x\in K\cap N$,
choose a basis $\vec{v}=\{v_1,\dots, v_k\}$ of $T_x(K\cap N)$. Extend
it to an oriented basis of $T_xK$ and of $T_xN$ as $\{v_1,\dots, v_k,
u_1,\dots, u_\ell\}$ and as $\{v_1,\dots,v_k,w_1,\dots,w_m\}$,
respectively. Then if $\{\vec{v}, \vec{w}, \vec{u}\}$ is the oriented
basis of $T_xM$, then $\vec{v}$ gives the orientation of $K\cap
N$. Otherwise, opposite orientation gives the orientation of $K\cap
N$.
\enddefinition

We come back to the homology of $\text{\rm SU}(n+1)$. For a sequence
$I$, let $\alpha_I\in H_*\left(\text{\rm SU}(n+1)\right)$ be the
Poincar\'e dual to $y_I$, that is $D(\alpha_I)=y_I$. Let $I^{c}$ be
the sequence complement to $I$ in $U$.  For two disjoint sequences
$I,J$ of strictly increasing integers, let $I\cup J$ be the combined
sequence of strictly increasing integers. Let
$\text{sgn}(I,J)\in\{\pm1\}$ be defined by
$$
x_{I\cup J}=\text{sgn}(I,J)x_I\cdot x_J
\tag2-2
$$
with respect to the Pontrjagin product. This is the same sign
resulting from converting the juxtaposition $IJ$ into a strictly
increasing sequence $I\cup J$. Note that
$\text{sgn}(I,J)=(-1)^{|I||J|}\text{sgn}(J,I)$. 

\proclaim{Lemma 2-3} \rom{(1)} For any sequence $I$, the homology class
$\alpha_I$ which is the Poincar\'e dual of $y_I$, is identified as 
$$
\alpha_I=(-1)^{\frac{|I|(|I|-1)}2}\text{\rm sgn}(I,I^c)x_{I^c}.
\tag2-3
$$

\rom{(2)} For two disjoint sequences $I,J$, the intersection product
is given by 
$$
\alpha_I\circ\alpha_J=\text{\rm sgn}(I,J)\alpha_{I\cup J}.
$$
Thus, if $I=(2i_1+1,2i_2+1,\cdots,2i_r+1)$ with $1\le
i_1<\cdots<i_r\le n$, we have 
$$
\alpha_I=\alpha_{2i_1+1}\circ\alpha_{2i_2+1}\circ\cdots\circ\alpha_{2i_r+1}. 
$$
\endproclaim
\demo{Proof} Since $D(\alpha_I)=y_I$, for any sequence $J$, we have 
$\langle y_J,\alpha_I\rangle=\langle y_J,D(\alpha_I)\cap x_U\rangle
=\langle y_J\cup y_I, x_U\rangle$, where $y_J\cup y_I=\text{\rm
sgn}(J,I)y_{I\cup J}$ and the last pairing is nonzero only when
$J=I^c$. Hence we have $\langle y_J,\alpha_I\rangle
=(-1)^{\frac{n(n-1)}2}\text{sgn}(I^c,I)\delta_{J,I^c}$, using
(2-1). This means that the homology class $\alpha_I$ is $x_{I^c}$ up
to a sign. Since $\langle
y_{I^c},x_{I^c}\rangle=(-1)^{\frac{|I^c|(|I^c|-1)}2}$, we see that
$\alpha_I=(-1)^{\frac{|I|(|I|-1)}2}\text{\rm sgn}(I,I^c)x_{I^c}$.

For the second formula, we compute:  
$$
D(\alpha_I\circ\alpha_J)=D(\alpha_I)\cup
D(\alpha_J)
=y_I\cup y_J=\text{\rm sgn}(I,J)y_{I\cup J}
=\text{\rm sgn}(I,J)D(\alpha_{I\cup J}). 
$$
This allows us to write $\alpha_I$ as an intersection product of
$\alpha_{2i+1}$'s. 
\qed
\enddemo

Since homology elements $\alpha_I$'s are dual to cohomology elements,
they are more suitable for computing intersection products than
$x_I$'s, and elements $x_I$'s are used for computing Pontrjagin
products. In homology, the intersection product is denoted by a circle
``$\circ$'' and the Pontrjagin product is denoted by a dot
``$\cdot$''.

Next we examine the relationship between the intersection product and
the Pontrjagin product. Properties of $\alpha_I$'s are best described
through Poincar\'e duality. 

\proclaim{Proposition 2-4} For any sequence $I$, we have 
$$
y_{2i+1}\cup D(x_{2i+1}\cdot\alpha_I)=
\cases 
y_I & \text{if }2i+1\in I,\\
0 & \text{if } 2i+1\notin I.
\endcases
$$
In other words, 
$$
x_{2i+1}\cdot\alpha_I=
\cases (-1)^{r-1}\alpha_{I-\{2i+1\}} & \text{if } 2i+1\in I,\\
0 & \text{if } 2i+1\notin I.
\endcases
$$
Here, when $2i+1$ appears in $I$, it appears as the $r$-th entry.
\endproclaim
\demo{Proof} Since $\alpha_I=\pm x_{I^c}$, if $2i+1\notin I$, then
$x_{2i+1}\cdot\alpha_I=0$. Suppose $2i+1\in I$. By (2-2) and (2-3),
we have 
$$
x_{2i+1}\cdot\alpha_I=(-1)^{\frac{|I|(|I|-1)}2}\text{\rm sgn}(I,I^c)
\text{\rm sgn}(\{2i+1\},I^c)x_{I^c\cup\{2i+1\}}. 
$$
On the other hand, by (2-3) we have 
$$
\alpha_{I-\{2i+1\}}=(-1)^{\frac{(|I|-1)(|I|-2)}2}
\text{\rm sgn}(I-\{2i+1\},I^c\cup\{2i+1\})x_{I^c}. 
$$
Comparing the signs, we get $x_{2i+1}\cdot\alpha_I=(-1)^{|I|-1}
\text{sgn}(I-\{2i+1\},\{2i+1\})\alpha_{I-\{2i+1\}}$. The last sign is
$(-1)^{|I|-r}$ if $2i+1$ appears in $I$ as the $r$-th entry. This
proves the second part. For the first part, we have 
$y_{2i+1}\cdot
D(x_{2i+1}\cdot\alpha_I)=y_{2i+1}\cdot(-1)^{r-1}y_{I-\{2i+1\}}=y_I$. 
\qed
\enddemo

This result means that the Pontrjagin product with $x_{2i+1}$
coincides, in the intersection ring 
$\Bbb H_*\left(\text{SU}(n+1)\right)
=\Lambda_{\Bbb Z}(\alpha_3,\cdots,\alpha_{2n+1})$, with the odd derivation
$\partial/\partial\alpha_{2i+1}$, where $\alpha_{2k+1}$'s are regarded
as odd variables to get the sign right. Thus for any sequence $I$, 
$$
x_{2i+1}\cdot \alpha_I
=\frac{\partial \alpha_I}{\partial \alpha_{2i+1}}, \quad 1\le i\le n. 
\tag2-4
$$
The next corollary is an immediate consequence of this fact. 

\proclaim{Corollary 2-5} For any two sequences $I,J$ not necessarily
disjoint, in $H_*\left(\text{\rm SU}(n+1)\right)$ we have 
$$
x_{2i+1}\cdot(\alpha_I\circ\alpha_J)=(x_{2i+1}\cdot\alpha_I)\circ\alpha_J
+(-1)^{|I|}\alpha_I\circ(x_{2i+1}\cdot\alpha_J).
$$
\endproclaim

\bigskip

\head
{\elevenbx \S3 Homology $S^1$-action on the Hopf algebra 
$H_*\left(L\text{\rm SU}(n+1)\right)$}
\endhead

\bigskip

For a Lie group $G$, its free loop space is topologically a product
$\xi: G\times \Omega G\cong LG$ in which $(g,\gamma)$ is mapped to $g\cdot
\gamma$ for an element $g\in G$ and a loop $\gamma$ based at the
identity element, using the multiplication $\mu: G\times G @>>> G$ in
$G$. Although the map $\xi$ commutes with diagonal maps, it does not
commute with group multiplications on both sides unless $G$ is
abelian.  Thus there may be a possibility that the Hopf algebra
structure of $H_*(LG)$ is different from the tensor product of Hopf
algebras $H_*(G)$ and $H_*(\Omega G)$. However, under mild assumptions
on $G$, we show that indeed the Hopf algebra structure $H_*(LG)$ is
isomorphic to the tensor product of two Hopf algebras $H_*(G)$ and
$H_*(\Omega G)$.

To compare the group structures of both sides of the map $\xi$, we
consider the conjugation map $\text{Ad}: \Omega G\times G @>>> \Omega
G$ given by $\text{Ad}(\gamma, g)=g^{-1}\cdot\gamma\cdot g$. Then the
following diagram commutes:
$$
\CD
(G\times \Omega G)\times (G\times \Omega G) @>{\xi\times \xi}>{\cong}>
LG\times LG \\
@VV{1\times1\times \phi\times 1}V   @VV{\mu}V  \\
G\times \Omega G\times G\times G\times \Omega G  @. LG \\
@VV{1\times\text{Ad}\times 1\times 1}V  @A{\cong}A{\xi}A  \\
G\times \Omega G\times G\times \Omega G @>{(\mu\times\mu)\circ(1\times
T\times 1)}>> G\times \Omega G 
\endCD
\tag3-1
$$
where $\phi$ is the diagonal map and $T$ is the switching map. To
study the homology of this diagram, we need to know the induced map
$\text{Ad}_*$ on homology, which turn out to be trivial in most
cases. 

\proclaim{Lemma 3-1} Assume that the Pontrjagin ring $H_*(G)$ is torsion
free and generated by odd degree elements. Then positive degree
elements of $H_*(G)$ acts trivially on $H_*(\Omega G)$ through
$\text{\rm Ad}_*$. Namely, for $e\in H_*(\Omega G)$ and $x\in H_*(G)$,  
$$
\text{\rm Ad}_*(e\otimes x)=
\cases
0  & \text{if } |x|>0,\\
e &  \text{if } x=1\in H_0(G).
\endcases
$$
\endproclaim
\demo{Proof} Since $H_*(G)\cong \Lambda_{\Bbb Z}(x_a\mid a\in A)$, 
where $A$ is an index set, we have $H_*(\Omega G)\cong \Bbb Z[e_a\mid
a\in A]$ with $|e_a|+1=|x_a|$. Then for each $a\in A$, the element
$\text{Ad}_*(e\times x_a)$ has an odd degree, and hence must vanish since
$\Bbb Z[e_a]$ has only even degree elements. Since $\text{Ad}_*$
defines a right action of the Pontrjagin ring $H_*(G)$ on $H_*(\Omega G)$,
we see that $\text{Ad}_*$ action vanishes for all positive degree
elements in $H_*(G)$.
\qed
\enddemo

An immediate consequence is that the composition
$(\text{Ad}\times 1)\circ(1\times \phi): \Omega G\times G @>>> \Omega
G\times G$ induces an identity map in homology. Thus, the commutative
diagram (3-1) implies the following result. 

\proclaim{Proposition 3-2} Assume that $H_*(G)$ is torsion free and
generated by odd degree elements. Then $H_*(LG)\cong H_*(G)\otimes
H_*(\Omega G)$ as Hopf algebras. 
\endproclaim
 
Next we examine the compatibility between the Hopf algebra structure of
$H_*(LG)$ and the homological $S^1$ action $\Delta: H_*(LG) @>>>
H_{*+1}(LG)$ given by $\Delta(z)=\Delta_*([S^1]\otimes z)$ for $z\in
H_*(LG)$. Here we consider $\Delta$ in the Pontrjagin ring $H_*(LG)$,
not in the loop homology algebra $\Bbb H_*(LG)$, which is the topic of
the next section.   

\proclaim{Proposition 3-3}  Let $G$ be a Lie group. 

\rom{(1)} The map $\Delta$ is a derivation on the Pontrjagin ring
$H_*(LG)$. Namely, for $z,w\in H_*(LG)$, 
$$
\Delta(z\cdot w)=\Delta (z)\cdot w+(-1)^{|z|}z\cdot\Delta(w).
$$

\rom{(2)} The map $\Delta$ commutes with the coalgebra map $\phi:
H_*(LG) @>>> H_*(LG)\otimes H_*(LG)$. Namely, if
$\phi(z)=\sum_iz_i'\otimes z_i''$ for $z\in H_*(LG)$, then 
$$
\phi\left(\Delta(z)\right)=\Delta\left(\phi(z)\right)
=\sum_i\left\{\Delta(z_i')\otimes z_i'' 
+(-1)^{|z_i'|}z_i'\otimes\Delta(z_i'')\right\}.
$$
\endproclaim
\demo{Proof} These are straightforward consequences of the following
commutative diagrams. 
$$
\CD
S^1\times (LG\times LG) @>{1\times\mu}>> S^1\times LG \\
@V{(1\times T\times 1)\circ(\phi\times 1\times 1)}VV   @V{\Delta}VV  \\
(S^1\times LG)\times (S^1\times LG)  @.       LG  \\
@V{\Delta\times\Delta}VV      @\vert  \\
LG\times LG @>{\mu}>> LG 
\endCD
\qquad
\CD
S^1\times LG @>{\Delta}>> LG \\
@V{(1\times T\times 1)\circ(\phi\times\phi)}VV    @V{\phi}VV \\
(S^1\times LG)\times (S^1\times LG) @>{\Delta\times \Delta}>> LG\times
LG
\endCD
$$

For a loop $\gamma\in LG$, let $\gamma_t$ denote the loop rotated by
$t\in S^1=\Bbb R/\Bbb Z$. Then the first diagram above says
$(\gamma\cdot\eta)_t=\gamma_t\cdot\eta_t$ for any two loops $\gamma$
and $\eta$. The second diagram says
$(\gamma,\gamma)_t=(\gamma_t,\gamma_t)$. 
\qed
\enddemo

If we can determine the effect of $\Delta$ on the generators of the
Pontrjagin ring $H_*(LG)$, then by derivation property of $\Delta$, we
can determine $\Delta$ on the entire ring $H_*(LG)$. 

To describe $\Delta$, we need differential operators $D_{2\ell}$
acting on the polynomial ring $\Bbb Z[e_2, e_4, \cdots, e_{2n}]$
introduced in (1-2). These differential operators have the following
properties. 

\proclaim{Proposition 3-4} Let $\phi$ be the coalgebra map in $\Bbb
Z[e_2, e_4, \cdots, e_{2n}]$ given in Theorem \rom{2-1}. 

\rom{(1)} The differential operators $D_2, D_4, \dots, D_{2n}$ are
characterized by 
$$
\phi(e^J)=1\otimes e^J+e_2\otimes D_2e^J+e_4\otimes
D_4e^J+\cdots+e_{2n}\otimes D_{2n}e^J+\!\!\sum_{I+K=J}\!\!a_{I,J}e^I\otimes e^K,
\tag3-2
$$
where in the last summation, $e^I$ is a product of at least two
elements. 

\rom{(2)} The operators $D_{2\ell}$'s mutually commute, that is, $[D_{2i},
D_{2j}]=0$ for $1\le i,j\le n$. 
\endproclaim
\demo{Proof} The formula in (1) can be checked directly. The
commutativity of $D_{2\ell}$'s can also be checked directly, but we
can also use the fact that the coalgebra map $\phi$ is cocommutative,
and we compare terms of the form $e_{2i}\otimes e_{2j}\otimes
(D_{2j}D_{2i}e^J)$ in the both sides of $(T\otimes 1) (1\otimes
\phi)\phi(e^J)=(1\otimes\phi)\phi(e^J)$. 
\qed
\enddemo 

To determine the effect of $\Delta$ on algebra generators of
$H_*\left(L\text{\rm SU}(n+1)\right)$, we use the coalgebra structure. 

\proclaim{Theorem 3-5} In the Pontrjagin ring 
$H_*\left(L\text{\rm SU}(n+1)\right)\cong
\Lambda(x_3,x_5,\dots,x_{2n+1})\otimes\Bbb Z[e_2,e_4, \dots, e_{2n}]$,
the effect of $\Delta$ on generators is given by 
$$
\Delta(x_{2\ell+1})=0, \qquad \Delta(e_{2\ell})=\sum_{i=1}^\ell
x_{2i+1}\cdot e_{2\ell-2i}, \qquad 1\le\ell\le n.
$$
The general formula of the action $\Delta$ is 
$$
\Delta(x_Ie^J)=\sum_{\ell=1}^n(x_{2\ell+1}\cdot x_I)\cdot(D_{2\ell}e^J).
\tag3-3
$$
\endproclaim
\demo{Proof} Since the map $x_{2k+1}: \Sigma\Bbb CP^k @>>>
\text{SU}(n+1) @>>> L\text{SU}(n+1)$ is a cycle consisting of constant
loops, the $S^1$-action on it is trivial and we have
$\Delta(x_{2k+1})=0$ for all $k$.  Since $x_{2k+1}$ is the adjoint of
the map $e_{2k}:\Bbb CP^k @>>> \Omega\text{SU}(n+1)$, $x_{2k+1}$ can
be thought of as the following composition:
$$
x_{2k+1}: S^1\times \Bbb CP^k @>{1\times e_{2k}}>> S^1\times
\Omega\text{SU}(n+1)
@>>> S^1\times L\text{SU}(n+1) @>{\Delta}>> L\text{SU}(n+1) @>{p}>>
\text{SU}(n+1),
$$
where $p$ is the projection to the base of loops. This means
$x_{2k+1}=p_*\left(\Delta_*([S^1]\otimes e_{2k})\right)
=p_*\left(\Delta(e_{2k})\right)$. Thus, in
$H_*\left(L\text{SU}(n+1)\right)$, we have
$\Delta(e_{2k})=x_{2k+1}+(\text{decomposables})$, where each term in
the decomposables has a nontrivial $e$-factor. We determine the
decomposable part by induction on $k$. When $k=1$, by dimensional
reason, we must have $\Delta(e_2)=x_3$. Now assume that
$\Delta(e_{2k})=\sum_{i=1}^kx_{2i+1}\cdot e_{2k-2i}$ for $1\le k\le
\ell-1$ for some $\ell\ge2$. Since $\Delta$ is a coalgebra map by
Proposition 3-3, we have
$\phi\left(\Delta(e_{2\ell})\right)=\Delta\left(\phi(e_{2\ell})\right)$. We compare
both sides. We let $\Delta(e_{2\ell})=x_{2\ell+1}+Y$, where $Y$ is a
decomposable element of homogeneous degree $2\ell+1$. Since
$\phi(e_{2\ell})=\sum_{k=0}^\ell e_{2k}\otimes e_{2\ell-2k}$ and 
$x_{2\ell+1}$ is primitive, the above identity gives the following
identity for $Y$:
$$
\phi(Y)=Y\otimes 1 + 1\otimes Y +
\sum_{k=1}^{\ell-1}\left(\Delta(e_{2k})\otimes
e_{2\ell-2k}+e_{2k}\otimes \Delta(e_{2\ell-2k})\right).
$$
Applying formula (3-2) to $Y$, and comparing it with the above
identity, we see that $Y$ must satisfy a system of PDEs given by
$D_{2k}(Y)=\Delta(e_{2\ell-2k})$ for $1\le k\le \ell-1$, and
$D_{2k}(Y)=0$ for $\ell\le k\le n$. Using the explicit expression of
operators $D_{2k}$ given in (1-2), the above system can be rewritten
in matrix form as 
$$
\left(
\matrix
1 & e_2 & \hdots & e_{2\ell-2} & \hdots & e_{2n-2} \\
  & 1   & \hdots & e_{2\ell-4} & \hdots & e_{2n-4} \\
  &     & \ddots & \vdots      & \ddots & \vdots  \\
  &     &        & 1           & \hdots & e_{2n-2\ell} \\
  &     &        &             & \ddots & \vdots      \\
  &     &        &             &        & 1
\endmatrix\right)
\left(
\matrix
\partial_2Y \\
\partial_4 Y \\
\vdots \\
\partial_{2\ell}Y \\
\vdots \\
\partial_{2n}Y
\endmatrix\right)
=\left(
\matrix 
\Delta(e_{2\ell-2}) \\
\Delta(e_{2\ell-4}) \\
\vdots \\
\Delta(e_2) \\
\vdots \\
0
\endmatrix
\right)
$$
where the coefficient matrix has zeros below the diagonal. Since the
coefficient matrix is invertible, we can uniquely determine partial
derivatives $\partial_{2k}Y=\partial Y/\partial e_{2k}$ for all $1\le
k\le n$. Since the constant term of $Y$ is zero, if the solution to
the system exists, it must be unique. On the other hand, it is easy to
check that $Y=\sum_{i=1}^{\ell-1}x_{2i+1}e_{2\ell-2i}$ satisfies the
system using the induction hypothesis. Thus, we have
$\Delta(e_{2\ell})=x_{2\ell+1}+Y=\sum_{i=1}^{\ell}x_{2i+1}e_{2\ell-2i}
=\sum_{i=1}^{\ell}x_{2i+1}D_i(e_{2\ell})$.  This proves the first
part.

The second part is a simple consequence of the derivation property of
$\Delta$ and differential operators $D_{2k}$, using the calculation
of $\Delta$ on algebra generators above. 
\qed
\enddemo

\bigskip

\head
{\elevenbx \S4 \quad Batalin-Vilkovisky Lie algebra structure on the loop
homology $\Bbb H_*\left(\text{LSU}(n+1)\right)$}
\endhead

\bigskip

In this section we describe the full BV algebra structure of the loop
homology of $\text{SU}(n+1)$. First we give a brief description of a 
general Batalin-Vilkovisky algebra. 

\definition{Definition} A Batalin-Vilkovisky algebra $A$ is a graded
commutative associative algebra with a degree $1$ operator $\Delta$
with $\Delta^2=0$ such that the bilinear product $\{\ ,\ \}$ on $A$
defined by 
$$
\Delta(ab)=\Delta(a)\cdot b +
(-1)^{|a|}a\cdot\Delta(b)+(-1)^{|a|}\{a,b\}, \quad a,b\in A
\tag4-1
$$
is a derivation in each variable in the sense that for every $a,b,c\in
A$, 
$$
\align
\{a,b\cdot c\}&=\{a,b\}\cdot c + (-1)^{|b|(|a|+1)}b\cdot\{a,c\},\\
\{a\cdot b,c\}&=a\cdot\{b,c\}+(-1)^{|b|(|c|+1)}\{a,c\}\cdot b.
\endalign
$$
\enddefinition 
It is a formal consequence that the bilinear product $\{\ ,\ \}$ defined above
satisfies identities for a Lie bracket of degree $1$, meaning 
$$
\align
\{a,b\}&=-(-1)^{(|a|+1)(|b|+1)}\{b,a\}, \\
\{a,\{b,c\}\}&=\{\{a,b\},c\}+(-1)^{(|a|+1)(|b|+1)}\{b,\{a,c\}\}.
\endalign
$$

We note that additively there is a bijective correspondence of elements between
the Pontrjagin ring $H_*\left(L\text{SU}(n+1)\right)$ and the loop
homology $\Bbb H_*\left(L\text{SU}(n+1)\right)$ preserving the
homological $S^1$-action $\Delta$, although the ring
structure and the grading are different between these two
rings. Specifically, this correspondence is given by (2-3).  Using
this bijection, we can translate the homological $S^1$-action $\Delta$
given in (3-3) in the Pontrjagin ring to the loop homology, and we get
the following formula: 
$$
\pmb\Delta(\alpha_Ie^J)=\sum_{\ell=1}^n\frac{\partial\alpha_I}
{\partial\alpha_{2\ell+1}}\cdot(D_{2\ell}e^J),
$$
for any sequences $I,J$, where $\partial/\partial \alpha_{2\ell+1}$ is
an odd derivation on odd variables $\alpha_{2k+1}$'s given in
(2.4). Here are some special cases of the above formula.
$$
\pmb\Delta(\alpha_I)=0, \qquad \pmb\Delta(e^J)=0, \qquad
\pmb\Delta(\alpha_{2\ell+1}e^J)=D_{2\ell}e^J, \qquad
\pmb\Delta(\alpha_{2\ell+1}e_{2\ell})=1, \quad 1\le\ell\le n. 
$$
We use the bold face $\pmb\Delta$ for the BV operator in the loop homology
to distinguish it from the $\Delta$ operator in the Pontrjagin ring. 

We can now compute the BV Lie bracket using the relation (4-1). For
arbitrary two elements $\alpha_Ie^J, \alpha_Ke^L\in \Bbb
H_*\left(L\text{SU}(n+1)\right)$, we have
$$
\aligned
\pmb\Delta(\alpha_Ie^J\circ\alpha_Ke^L)
&=\sum_{i=1}^n\frac{\partial}{\partial
\alpha_{2i+1}}(\alpha_I\alpha_K)\cdot D_{2i}(e^Je^L) \\
&=\sum_{i=1}^n\left(\frac{\partial \alpha_I}
{\partial\alpha_{2i+1}}\cdot D_{2i}e^J\right)\cdot\alpha_Ke^L
+(-1)^{|I|}\alpha_Ie^J\cdot
\sum_{i=1}^n\left(\frac{\partial \alpha_K}{\partial
\alpha_{2i+1}}\cdot D_{2i}e^L\right) \\
&\qquad +\sum_{i=1}^n\left(\frac{\partial\alpha_I}{\partial
\alpha_{2i+1}}\cdot D_{2i}e^L\right)\circ \alpha_Ke^J
+(-1)^{|I|}\alpha_Ie^L\circ\sum_{i=1}^n\left(\frac{\partial\alpha_K}{\partial
\alpha_{2i+1}}\cdot D_{2i}e^J\right)\\
=\pmb\Delta(\alpha_Ie^J)\circ&\alpha_Ke^L
+(-1)^{|I|}\alpha_Ie^J\circ\pmb\Delta(\alpha_Ke^L)
+\pmb\Delta(\alpha_Ie^L)\circ\alpha_Ke^J
+(-1)^{|I|}\alpha_Ke^L\circ\pmb\Delta(\alpha_Ke^J).
\endaligned
$$
Comparing with the formula (4-1), we see that the last two terms give
the BV Lie bracket. Collecting our results for the loop homology, we get
the next Theorem.

\proclaim{Theorem 4-1} In the loop homology $\Bbb H_*\left(L\text{\rm
SU}(n+1)\right)
\cong\Lambda_{\Bbb Z}(\alpha_3,\alpha_5,\cdots,\alpha_{2n+1})\otimes
\Bbb Z[e_2,e_4, \cdots, 2_{2n}]$, the BV operator $\pmb\Delta$ and BV
Lie bracket are given by the following formulae\rom{:}
$$
\gathered
\pmb\Delta(\alpha_Ie^J)=\sum_{\ell=1}^n\frac{\partial\alpha_I}
{\partial\alpha_{2\ell+1}}\cdot(D_{2\ell}e^J), \\
(-1)^{|I|}\left\{\alpha_Ie^J,\alpha_Ke^L\right\}
=\pmb\Delta(\alpha_Ie^L)\circ\alpha_Ke^J
+(-1)^{|I|}\alpha_Ie^L\circ\pmb\Delta(\alpha_Ke^J).
\endgathered
$$
\endproclaim

It is rather surprising that the BV Lie bracket in our loop homology
is given precisely by the $\pmb\Delta$-derivation formula with
$e$-terms switched. We want to understand what is behind this rather
intriguing formula. 

So far, we have been examining the Batalin-Vilkovisky algebra
structure in loop homology, but loop homology has more structures. For
example, it admits so called string operations. To understand the BV
bracket, a relevant structure in loop homology not contained in BV
algebra structure is the {\it coalgebra} structure in
$H_*\left(\Omega\text{SU}(n+1)\right)$. Especially, its primitive
elements play an essential role. In Theorem 2-1, we described its coalgebra
structure. Primitive elements are spanned by $s$-classes
$s_{2\ell}(e_2,e_4,\cdots, e_{2\ell})$ for $1\le \ell\le n$, defined
inductively by Newton's formula:
$$
s_{2\ell}-e_2s_{2\ell-2}+e_4s_{2\ell-4}-\cdots
+(-1)^{\ell-1}e_{2\ell-2}s_2+(-1)^{\ell}\ell e_{2\ell}=0, \quad
1\le\ell\le n. 
$$
Note that $s_{2\ell}=(-1)^{\ell-1}\ell
e_{2\ell}+(\text{decomposables})$. We extend the ground ring from
$\Bbb Z$ to $\Bbb Q$ and we define elements 
$$
h_{2\ell}=(-1)^{\ell-1}\frac{s_{2\ell}}{\ell}=e_{2\ell}+(\text{decomposables})
\in\Bbb Q[e_2,\cdots,e_{2n}], \quad 1\le\ell\le n.
$$
Then the algebra $\Bbb Q[e_2,e_4,\cdots,e_{2n}]=\Bbb
Q[h_2,h_4,\cdots,h_{2n}]$ is primitively generated by $h_{2\ell}$'s
and we may regard $\{h_2,h_4,\dots,h_{2n}\}$ as new coordinate
functions on $n$-dimensional vector space $\Bbb Q^n$. As such,
differential operators $\partial/\partial h_{2\ell}$ makes sense for
functions in $e_2,\cdots ,e_{2n}$. 

\proclaim{Proposition 4-2} For each $1\le\ell\le n$, the derivation
$\partial/\partial h_{2\ell}$ preserves the integral subring $\Bbb Z[e_2,e_4,
\cdots, e_{2n}]$, and is given by 
$$
\frac{\partial}{\partial h_{2\ell}}=D_{2\ell}
=\sum_{k=\ell}^ne_{2k-2\ell}\frac{\partial}{\partial e_{2k}}. 
$$
Furthermore, as a function in $e_2, \cdots, e_{2n}$, all the partial
derivatives of $h_{2\ell}$ have integral coefficients:
$$
\frac{\partial h_{2\ell}}{\partial e_{2k}}\in\Bbb
Z[e_2,e_4,\cdots,e_{2n}], \quad 1\le k,\ell\le n. 
$$
\endproclaim
\demo{Proof} For the first part, we need to show that the functions
$h_{2\ell}(e_2,\cdots,e_{2\ell})$ satisfy the system
$D_{2k}h_{2\ell}=\delta_{k,\ell}$ for all $1\le k,\ell\le n$. It is
possible to verify this directly by induction, but it is rather long
and tedious. Here we take the following approach, which gives a quick
proof.  Recall that the differential operators $D_{2k}$ are defined
in terms of coalgebra map $\phi$ in $\Bbb Z[e_2,\cdots,2_{2n}]$ as in
(3-2). Extending $\phi$ to over $\Bbb Q$, we have
$$
\phi(h_{2\ell})=1\otimes h_{2\ell} + e_2\otimes D_2(h_{2\ell}) +
e_4\otimes D_4(h_{2\ell})+\cdots+e_{2n}\otimes D_{2n}(h_{2\ell}) 
+\sum_{I,J}a_{I,J}e^I\otimes e^J.
$$
On the other hand, the element $h_{2\ell}$ is primitive. So we have 
$\phi(h_{2\ell})=1\otimes h_{2\ell}+h_{2\ell}\otimes 1
=1\otimes h_{2\ell}+e_{2\ell}\otimes 1 + (\text{decomposables})\otimes
1$. Comparing these two identities, we see that
$D_{2k}(h_{2\ell})=\delta_{k,\ell}$ for all $k,\ell$.   

For the second part, writing the system
$D_{2k}h_{2\ell}=\delta_{k,\ell}$ for a given $\ell$ in a matrix form,
we get 
$$
\left(
\matrix
1 & e_2 & \hdots & e_{2\ell-2} & \hdots & e_{2n-2} \\
  & 1   & \hdots & e_{2\ell-4} & \hdots & e_{2n-4} \\
  &     & \ddots & \vdots      & \ddots & \vdots  \\
  &     &        & 1           & \hdots & e_{2n-2\ell} \\
  &     &        &             & \ddots & \vdots      \\
  &     &        &             &        & 1
\endmatrix\right)
\left(
\matrix
\partial_2h_{2\ell} \\
\partial_4 h_{2\ell} \\
\vdots \\
\partial_{2\ell}h_{2\ell} \\
\vdots \\
\partial_{2n}h_{2\ell}
\endmatrix\right)
=\left(
\matrix
0 \\
0 \\
\vdots \\
1 \\
\vdots \\
0
\endmatrix
\right)
$$
where in the right hand side, $1$ is in the $\ell$th position. Since
the coefficient matrix is invertible in the integral ring $\Bbb
Z[e_2,\cdots,e_{2n}]$, we see that all the partial derivatives of
$h_{2\ell}$ are in this integral ring.
\qed
\enddemo

To understand BV Lie bracket in the loop homology, we take the point
of view that $\Bbb H_*^{\Bbb Q}=\Bbb
H_*\left(\text{SU}(n+1)\right)\otimes\Bbb Q$ is the ring of functions
on the super vector space $\Bbb Q^{n|n}$ of dimension $n|n$ in which
$h_2, h_4,\cdots, h_{2n}$ are even coordinate functions, and $\alpha_{3},
\alpha_5, \cdots, \alpha_{2n+1}$ are odd coordinate functions. We can do usual
differential geometry in this super setting and we can consider vector
fields and differentials. In particular, we can consider an odd
symplectic form
$$
\omega=\sum_{\ell=1}^nd\alpha_{2\ell+1}\wedge dh_{2\ell}
\in\bigwedge\!\!{}^2\Bigl(\bigoplus_{\ell=1}^n\Bbb H_*^{\Bbb Q}
d\alpha_{2\ell+1}\oplus
\bigoplus_{\ell=1}^n\Bbb H_*^{\Bbb Q}dh_{2\ell}\Bigr).
$$
The associated Poisson bracket $[\ \ ,\ \ ]_{\omega}: \Bbb H_*^{\Bbb
Q}\otimes\Bbb H_*^{\Bbb Q} @>>> \Bbb H_*^{\Bbb Q}$, is given by
$[F,G]=X_FG$ for $F,G\in\Bbb H_*^{\Bbb Q}$, where the vector field
$X=X_F$ is defined by the relation $\iota_{X}\omega=dF$. We follow
\cite{4} for a sign convention. For example, the exterior derivative
is given by
$d=\sum_{\ell}d\alpha_{2\ell+1}\otimes(\partial/\partial\alpha_{2\ell+1})
+\sum_{\ell}dh_{2\ell}\otimes(\partial/\partial h_{2\ell})$, and for
$F=\alpha_Ih^J\in\Bbb H_*^{\Bbb Q}$ the associated vector field $X_F$
is given by
$$
X_F=\sum_{\ell=1}^n\left(
(-1)^{|I|}\frac{\partial F}{\partial
\alpha_{2\ell+1}}\cdot\frac{\partial}{\partial h_{2\ell}}
+\frac{\partial F}{\partial h_{2\ell}}\cdot
\frac{\partial}{\partial\alpha_{2\ell+1}}\right).
$$

\proclaim{Theorem 4-3} \rom{(1)} The differentials $dh_{2\ell}$ and
derivations $\partial/\partial h_{2\ell}$ for $1\le\ell\le n$ are
defined over $\Bbb Z$. Consequently, the symplectic form and the
Poisson bracket are defined over $\Bbb Z$\rom{:}
$$
\omega\in \bigwedge\!\!{}^2\Bigl(\bigoplus_{\ell=1}^n\Bbb H_*d\alpha_{2\ell+1}\oplus
\bigoplus_{\ell=1}^n\Bbb H_*dh_{2\ell}\Bigr), \quad\text{and}\quad
[\ \ ,\ \ ]_{\omega}: \Bbb H_*\otimes\Bbb H_* @>>> \Bbb H_*.
$$

\rom{(2)} The Poisson bracket $[\ \ ,\ \ ]_{\omega}$ defined on $(\Bbb
H_*,\omega)$ coincides with the Batalin-Vilkovisky Lie
bracket. Namely, 
$$
\left\{\alpha_Ie^J,\alpha_Ke^L\right\}
=\left[\alpha_Ie^J,\alpha_Ke^L\right]_{\omega},
$$
for $\alpha_Ie^J,\alpha_Ke^L\in \Bbb H_*$. 
\endproclaim
\demo{Proof} The first part is a consequence of Proposition 4-2. For
the second part, using $\partial/\partial h_{2\ell}=D_{2\ell}$, we have 
$$
\aligned
[\alpha_Ie^J,\alpha_Ke^L]_{\omega}&=X_{\alpha_Ie^J}(\alpha_Ke^L) \\
&=\sum_{\ell=1}^n
\left((-1)^{|I|}\frac{\partial \alpha_I}{\partial
\alpha_{2\ell+1}}\cdot\frac{\partial e^L}{\partial h_{2\ell}}\cdot \alpha_Ke^J
+\alpha_Ie^L\cdot \frac{\partial \alpha_K}{\partial\alpha_{2\ell+1}}
\frac{\partial e^J}{\partial h_{2\ell}}\right) \\
&=(-1)^{|I|}\left(\sum_{\ell=1}^n
\frac{\partial \alpha_I}{\partial\alpha_{2\ell+1}}\cdot
D_{2\ell}e^L\right) \cdot \alpha_Ke^J +\alpha_Ie^L\cdot
\left(\sum_{\ell=1}^n\frac{\partial \alpha_K}{\partial\alpha_{2\ell+1}}\cdot
D_{2\ell}e^J\right) \\
&=(-1)^{|I|}\pmb\Delta(\alpha_Ie^L)\circ\alpha_Ke^J
+\alpha_Ie^L\circ\pmb\Delta(\alpha_Ke^J) \\
&=\left\{\alpha_Ie^J,\alpha_Ke^L\right\}
\endaligned
$$
This proves the second part.
\qed
\enddemo

Let $\text{Ham}(\Bbb H_*^{\Bbb Q})=\{x\in\text{Der}(\Bbb H_*^{\Bbb
Q})\mid \Cal L_X\omega=0\}$ be the set of all Hamiltonian derivations,
where $\Cal L_X$ is the Lie derivative by $X$.  Since
$d[F,G]=\iota_{[X_{F},X_G]}\omega$ for $F,G\in\Bbb H_*^{\Bbb Q}$,
the correspondence $\Bbb H_*^{\Bbb Q} @>>> \text{Ham}(\Bbb H_*^{\Bbb
Q})$ given by $F \mapsto X_F$ is a super Lie algebra homomorphism with
kernel $\Bbb Q\cdot 1$. Since our symplectic manifold $\Bbb
Q^{n|n}$ is affine, the above correspondence is onto. This gives a
description of $\Bbb H_*^{\Bbb Q}$ modulo $\Bbb Q\cdot 1$ as the super
Lie algebra of Hamiltonian derivations. 

Now consider a vector subspace $Q$ of $\Bbb H_*^{\Bbb Q}$ spanned by
quadratic elements $\alpha_{2\ell+1}h_{2k},
\alpha_{2\ell+1}\alpha_{2j+1}$, and  $h_{2\ell}h_{2k}$ with $1\le \ell,
k,j\le n$, $\ell\not=j$. The corresponding Hamiltonian derivations are 
$$
\aligned
X_{\alpha_{2\ell+1}h_{2k}}
&=\alpha_{2\ell+1}\frac{\partial}{\partial
\alpha_{2k+1}}-h_{2k}\frac{\partial}{\partial h_{2\ell}},\\
X_{\alpha_{2\ell+1}\alpha_{2j+1}}
&=\alpha_{2k+1}\frac{\partial}{\partial h_{2\ell}}
-\alpha_{2\ell+1}\frac{\partial}{\partial h_{2k}}, \\
X_{h_{2\ell}h_{2k}}
&=h_{2\ell}\frac{\partial}{\partial \alpha_{2k+1}}
+h_{2k}\frac{\partial}{\partial \alpha_{2\ell+1}}.
\endaligned
$$
Let $V=\bigoplus_{\ell=1}^n\Bbb Q\,\alpha_{2\ell+1}\oplus
\bigoplus_{\ell=1}^n\Bbb Q\,h_{2\ell}$. The above
derivations induce linear maps $f: V @>>> V$ annihilating the odd
symplectic form $\omega$ in the sense that
$$
\sum_{\ell=1}^n\left\{d\left(f(\alpha_{2\ell+1})\right)\wedge
dh_{2\ell}+d\alpha_{2\ell+1}\wedge d\left(f(h_{2\ell})\right)\right\}=0, 
$$
and by counting the dimension, they form a basis for a super Lie
algebra $\pi sp(n|n)$. This proves Corollary D in the introduction.
 
We consider the Lie algebra structure of the loop homology of
$\text{SU}(n+1)$ from another point of view. For $0\le p\le n$, let
$L(p)$ be the subgroup of the (integral) loop homology $\Bbb H_*=\Bbb
H_*\left(L\text{SU}(n+1)\right)$ spanned by elements $\alpha_Ie^J$
with $|I|=p$ and arbitrary $J$. Then $\Bbb H_*$ is a direct sum $\Bbb
H_*=\bigoplus_{p=1}^n L(p)$. In view of Theorem 4-1, BV Lie bracket
decreases the number of $\alpha$-factors by $1$, and we have
$$
\{L(p),L(q)\}\subset L(p+q-1), \quad 0\le p,q\le n,
$$
where $L(p)=0$ when $p>n$. These subgroups have the following
properties:
\roster
\item $L(1)$ is a Lie subalgebra.
\item $L(0)$ is an abelian Lie subalgebra. 
\item $L(p)$ is a module over the Lie algebra $L(1)$ for $0\le p\le n$. 
\item $L(1)$ is the Lie algebra of all derivations on $L(0)=\Bbb
Z[e_2,e_4,\cdots,e_{2n}]$. 
\endroster

Only the last statement needs to be verified. By Theorem 4-1, the
action of $L(1)$ on $L(0)$ is given by
$\{\alpha_{2\ell+1}e^J,e^L\}=-e^J\cdot D_{2\ell}e^L$. Thus the effect
of $\alpha_{2\ell+1}e^J$ on $L(0)$ is the same as the derivation
-$e^JD_{2\ell}$. Since over $\Bbb Z[e_2,\dots, e_{2n}]$, operators
$D_2,\cdots, D_{2n}$ span the same set of derivations as operators
$\partial/\partial e_2, \cdots, \partial/\partial e_{2n}$ do, $L(1)$
is the Lie algebra of all derivations on $L(0)$.

\bigskip

\head
{\elevenbx \S5 Batalin-Vilkovisky algebra structure on the loop homology of
complex Stiefel manifolds}
\endhead

\bigskip

Let $V_k(\Bbb C^n)$ be the complex Stiefel manifold of complex
orthonormal $k$-frames in $\Bbb C^n$ with respect to the standard
Hermitian pairing in $\Bbb C^n$. We determine the
Batalin-Vilkovisky algebra structure of its loop homology $\Bbb
H_*\left(LV_k(\Bbb C^n)\right)$. It is well known that the complex
Stiefel manifolds are homogeneous spaces: 
$$
V_{n+1-k}(\Bbb C^{n+1})
=\text{SU}(n+1)/\text{SU}(k), \quad\text{ for }1\le k\le n, \quad
V_{n+1}(\Bbb C^{n+1})=\text{U}(n+1). 
$$
Note that $V_n(\Bbb C^{n+1})=\text{SU}(n+1)$, and $V_1(\Bbb
C^{n+1})=S^{2n+1}$. Since $\text{U}(n+1)=\text{SU}(n+1)\times S^1$ as
topological spaces, we have an isomorphism as BV algebras:
$$
\Bbb H_*\left(L\text{U}(n+1)\right)\cong
\Bbb H_*\left(L\text{SU}(n+1)\right)\otimes\Bbb H_*\left(LS^1\right)
$$ 
Thus for the remainder of this section, we consider the BV algebra
$\Bbb H_*\left(LV_{n+1-k}(\Bbb C^{n+1})\right)$ with $1\le k\le n$. We
will see that its description is very similar to the one for $\Bbb
H_*\left(L\text{SU}(n+1)\right)$. 

First we review the homology and cohomology of the homogeneous space
$\text{SU}(n+1)/\text{SU}(k)$. We define a set $\Cal S_k$ of sequences
by 
$$
\Cal S_k=\{(2i_1+1,\cdots,2i_r+1)\mid k\le i_1<\cdots <i_r\le n\}. 
$$
Let $U=U_k=(2k+1,2k+3,\cdots,2n+1)\in\Cal S_k$. Recall that for $1\le
i\le n$ we have cycles $x_{2i+1}:\Sigma\Bbb CP^i @>>>
\text{SU}(n+1)$ for $\text{SU}(n+1)$. For $I\in\Cal S_k$, we
define a cycle $x_I$ in $\text{SU}(n+1)/\text{SU}(k)$ by 
$$
x_I=x_{2i_1+1}\times \cdots \times x_{2i_r+1}:
\prod_{t=1}^r\Sigma\Bbb CP^{i_t} @>>> \text{SU}(n+1) @>>>
\text{SU}(n+1)/\text{SU}(k).
$$
In view of the cell decomposition of $\text{SU}(n+1)/\text{SU}(k)$ using
the product of top cells of $x_{2i+1}$'s, the above cycles span the
entire homology group. Thus
$$
H_*\left(\text{SU}(n+1)/\text{SU}(k)\right)\cong\bigoplus_{I\in\Cal
S_k}\Bbb Zx_I.
$$
The element $x_{U}$ defines the orientation class.  Since
$\text{SU}(n+1)/\text{SU}(k)$ is not an $H$-space, there is no
Pontrjagin ring structure in this homology group. However, there is a
coalgebra structure given by $\phi(x_I)=\sum_{J\cup
K=I}\text{sgn}(J,K)x_J\otimes x_K$. Primitive elements are spanned by
$x_{2\ell+1}$ for $k\le\ell\le n$. Let $y_{2\ell+1}$ be a cohomology
element dual to the primitive $x_{2\ell+1}$ with respect to the basis
$\{x_I\}$. Then it is straightforward to see that
$$
H^*\left(\text{SU}(n+1)/\text{SU}(k)\right)\cong
\Lambda_{\Bbb Z}(y_{2k+1},y_{2k+3},\cdots,y_{2n+1}).
$$
Let $D:H_*\left(\text{SU}(n+1)/\text{SU}(k)\right) @>{\cong}>>
H^{d-*}\left(\text{SU}(n+1)/\text{SU}(k)\right)$ be the Poincar\'e duality
map with respect to the orientation class $x_U$ so that that
$D(x_I)\cap x_U=x_I$ for all $I\in\Cal S_k$, where $d$ is the
dimension of the homogeneous space. As before, let $\alpha_I$
be the homology class Poincar\'e dual to $y_I$ for $I\in\Cal
S_k$. Then as for the case of $\text{SU}(n+1)$, we have 
$$
\alpha_I=(-1)^{\frac{|I|(|I|-1)}2}\text{sgn}(I,I^c)x_{I^c},
\tag5-1
$$
where the complement $I^c$ is taken in $U=U_k$.  The intersection ring
is then given by 
$$
\Bbb H_*\left(\text{SU}(n+1)/\text{SU}(k)\right)
\cong\Lambda_{\Bbb Z}(\alpha_{2k+1},
\alpha_{2k+3},\cdots,\alpha_{2n+1}),
$$
where the generator $\alpha_{2\ell+1}=(-1)^{\ell-k}x_{U-\{2\ell+1\}}$
has codimension $2\ell+1$ in $\text{SU}(n+1)/\text{SU}(k)$, and has
degree $-2\ell-1$ in the intersection ring. 

Next we describe the homology of the based loop space 
$\Omega\left(\text{SU}(n+1)/\text{SU}(k)\right)$. Using the map
$e_{2\ell}$ for the $\Omega\text{SU}(n+1)$, we let $e_{2\ell}: \Bbb
CP^{\ell} @>{e_{2\ell}}>> \Omega\text{SU}(n+1) @>>>
\Omega\left(\text{SU}(n+1)/\text{SU}(k)\right)$, using the same
symbol. Then $x_{2\ell+1}:\Sigma\Bbb CP^{\ell} @>>>
\text{SU}(n+1)/\text{SU}(k)$ is adjoint to $e_{2\ell}:\Bbb CP^{\ell}
@>>> \text{SU}(n+1)/\text{SU}(k)$. By a standard argument, we have 
$$
H_*\left(\Omega\left(\text{SU}(n+1)/\text{SU}(k)\right)\right)\cong
\Bbb Z[e_{2k}, e_{2k+2}, \cdots, e_{2n}],
$$
and the map $p_*:H_*\left(\Omega(\text{SU}(n+1))\right) @>>>
H_*\left(\Omega\left(\text{SU}(n+1)/\text{SU}(k)\right)\right)$ is a
quotienting map by an ideal $(e_2,e_4,\cdots,e_{2k-2})$. Since the
above map $p_*$ is a Hopf algebra map, all the Hopf algebra structure
of $H_*\left(\Omega\left(\text{SU}(n+1)/\text{SU}(k)\right)\right)$ is
obtained by setting $e_2=\cdots=e_{2k-2}=0$ in the structure maps of
$H_*\left(\Omega\left(\text{SU}(n+1)\right)\right)$.

Next, we consider the fibre bundle
$$
\Omega\left(\text{SU}(n+1)/\text{SU}(k)\right)
 @>>> L\left(\text{SU}(n+1)/\text{SU}(k)\right)
@>>> \text{SU}(n+1)/\text{SU}(k). 
$$
The associated homology and loop homology spectral sequences
$$
\aligned
E_{*,*}^2&=H_*\left(\text{SU}(n+1)/\text{SU}(k)\right)
\otimes H_*\left(\Omega\left(\text{SU}(n+1)/\text{SU}(k)\right)\right) 
\Longrightarrow
H_*\left(L\left(\text{SU}(n+1)/\text{SU}(k)\right)\right), \\
E_{*,*}^2&=\Bbb H_*\left(\text{SU}(n+1)/\text{SU}(k)\right)
\otimes H_*\left(\Omega\left(\text{SU}(n+1)/\text{SU}(k)\right)\right)
\Longrightarrow
\Bbb H_*\left(L\left(\text{SU}(n+1)/\text{SU}(k)\right)\right),
\endaligned
$$
collapse by a standard argument, although the fibre bundle itself is
not a product. Thus 
$$
H_*\left(L\left(\text{SU}(n+1)/\text{SU}(k)\right)\right)\cong
\bigoplus_{I,J}\Bbb Zx_Ie^J,
$$
where $I\in\Cal S_k$ and $J$ runs over all sequences
$J=(j_k,\cdots,j_n)$ with nonnegative integer entries. This homology
is only a coalgebra and does not have a ring structure. We can think
of the element $x_Ie^J$ as the image of the element of the same name
in $H_*\left(L\text{SU}(n+1)\right)$ under the projection map $p_*$. On
the other hand, the loop homology algebra is given by
$$
\Bbb H_*\left(L\left(\text{SU}(n+1)/\text{SU}(k)\right)\right)\cong
\Lambda_{\Bbb Z}(\alpha_{2k+1},\cdots,\alpha_{2n+1})\otimes\Bbb
Z[e_{2k},\cdots,2_{2n}],
\tag5-2
$$ 
but there is no coalgebra structure. Note that the map 
$p:\text{SU}(n+1) @>>> \text{SU}(n+1)/\text{SU}(k)$ does not induce a
homomorphism in loop homology.  

To determine the BV operator $\pmb\Delta$ in loop homology, we first
determine homological $S^1$-action in ordinary homology, and then
translate it into loop homology using (5-1). Since $S^1$ action
commutes with the projection map $p: L\text{SU}(n+1) @>>> 
L\left(\text{SU}(n+1)/\text{SU}(k)\right)$, the homological $S^1$
action in $H_*\left(L\left(\text{SU}(n+1)/\text{SU}(k)\right)\right)$
is obtained by projection from the formula of $\Delta$ in
$H_*\left(\text{SU}(n+1)\right)$. Thus for an element 
$x_Ie^J\in H_*\left(L\left(\text{SU}(n+1)/\text{SU}(k)\right)\right)$,
we lift it to an element $x_Ie^J$ in $H_*\left(\text{SU}(n+1)\right)$,
apply $\Delta$, then project down by $p_*$:
$$
\Delta(x_Ie^J)=p_*\left(\Delta(x_Ie^J)\right)
=\sum_{\ell=k}^np_*(x_{2\ell+1}\cdot x_I)\cdot D_{2\ell}^{(k)}e^J,
$$
where $D_{2\ell}^{(k)}$ is the differential operator obtained from
$D_{2\ell}$ by setting $e_2=\cdots=e_{2k-2}=0$. We can rewrite this
formula in terms of elements $\alpha_{2\ell+1}$'s. As in (2-4) for
$\text{SU}(n+1)$ case, we can show $p_*(x_{2\ell+1}\cdot
\alpha_I)={\partial
\alpha_I}/\partial\alpha_{2\ell+1}$.  Thus the BV structure of the
loop homology of the complex Stiefel manifold $V_{n+1-k}(\Bbb
C^{n+1})=\text{SU}(n+1)/\text{SU}(k)$ is formally entirely analogous
to the BV structure of $\Bbb H_*\left(L\text{SU}(n+1)\right)$, and we
obtain the following statement. 

\proclaim{Theorem 5-1} In the loop homology 
$\Bbb H_*\left(L\left(\text{\rm SU}(n+1)/\text{\rm SU}(k)\right)\right)$ of
the complex Stiefel manifold given in \rom{(5-2)}, the BV operator
$\pmb\Delta$ is given by 
$$
\pmb\Delta(\alpha_Ie^J)
=\sum_{\ell=k}^n\frac{\partial \alpha_I}{\partial\alpha_{2\ell+1}}
\cdot D_{2\ell}^{(k)}e^J,
$$
where differential operators $D_{2\ell}^{(k)}$ on the polynomial ring
$\Bbb Z[e_{2k},\cdots, e_{2n}]$ is given by 
$$
D_{2\ell}^{(k)}=
\cases 
\dfrac{\partial}{\partial e_{2\ell}} +
\dsize\sum\limits_{j=k}^{n-\ell}e_{2j}\dfrac{\partial}{\partial e_{2\ell+2j}}, &
\text{if } \ell\le n-k,\\
\dfrac{\partial}{\partial e_{2\ell}} & 
\text{if } \ell> n-k. 
\endcases
$$

The BV Lie bracket is given by 
$$
\aligned
\{\alpha_Ie^J,\alpha_Ke^L\}
&=(-1)^{|I|}\pmb\Delta(\alpha_Ie^L)\circ\alpha_Ke^J
+\alpha_Ke^L\circ\pmb\Delta(\alpha_Ke^J) \\
&=[\alpha_Ie^J,\alpha_Ke^L]_{\omega},
\endaligned
$$
where $\omega=\sum_{\ell=k}^nd\alpha_{2\ell+1}\wedge
dh_{2\ell}^{(k)}$ is the odd symplectic form associated with primitive
elements $h_{2\ell}^{(k)}$ for $k\le\ell\le n$ in $\Bbb
Q[e_{2k},\cdots,e_{2n}]$, and $[\cdot ,\cdot ]_{\omega}$ is the Poisson
bracket associated to $\omega$. 
\endproclaim

The proof is entirely analogous to the $\text{SU}(n+1)$ case. 

As a final result, we specialize to the case $V_1(\Bbb
C^{n+1})=\text{SU}(n+1)/\text{SU}(n)=S^{2n+1}$. 

\proclaim{Theorem 5-2} In the loop homology of an odd dimensional
sphere
$$
\Bbb H_*(LS^{2n+1})\cong\Lambda(\alpha)\otimes\Bbb Z[h],
$$ 
where $|\alpha|=-2n-1$ and $|h|=2n$, the BV operator is given by 
$$
\pmb\Delta(h^{\ell})=0, \qquad 
\pmb\Delta(\alpha h^{\ell})=\ell h^{\ell-1}, \quad
\ell\ge0.
$$
As a BV Lie algebra, $\Bbb H_{\text{\rm even}}=\Bbb Z[h]$ is an
abelian Lie subalgebra, and $\Bbb H_{\text{\rm odd}}=\alpha\Bbb
Z[h]$ is the Lie algebra of derivations on the polynomial ring
$\Bbb Z[h]$. Namely, 
$$
\{\alpha h^{k+1},\alpha h^{\ell+1}\}=(k-\ell)\alpha h^{k+\ell+1},\quad
\{\alpha h^{k+1},h^m\}=-mh^{k+m},\quad k,\ell\ge-1, m\ge0.
$$
\endproclaim

\definition{Remark}
If we let $L(k)=\alpha h^{k+1}$ for $k\ge-1$, then these
elements satisfy 
$$
\{L(k),L(\ell)\}=(k-\ell)L(k+\ell), \quad k,\ell\ge-1.
$$
This is the commutation relation of the positive half of the Virasoro
algebra. Its generalization is the Lie algebra $L(1)$ acting as the
Lie algebra of derivations on the polynomial algebra $\Bbb Z[e_2,e_4,
\cdots, e_{2n}]$ discussed at the end of section 4.
\enddefinition

The BV algebra in Theorem 5-2 is a building block of the loop homology
of $\text{SU}(n+1)$, as we saw in Corollary C.

\enddocument